\documentclass[10pt]{amsart}
\usepackage[utf8]{inputenc}

\usepackage{amssymb,amsthm,amsmath,mathrsfs}
\usepackage{enumerate}
\usepackage{graphicx,color}
\usepackage{setspace}
\usepackage[hidelinks]{hyperref}
\usepackage{float}

\newcommand{\dd}{\mathrm{d}}
\newcommand{\E}{\mathbb{E}}
\newcommand{\1}{\textbf{1}}
\newcommand{\R}{\mathbb{R}}

\newcommand{\p}[1]{\mathbb{P}\left( #1 \right)}

\DeclareMathOperator{\vol}{vol}

\usepackage{xpatch}
\makeatletter
\def\thm@space@setup{%
  \thm@preskip=12pt plus 0pt minus 0pt
  \thm@postskip=0pt plus 0pt minus 0pt
}
\xpatchcmd{\proof}{6\p@\@plus6\p@\relax}{\z@skip}{}{}
\makeatother

\newtheorem{theorem}{Theorem}
\newtheorem{lemma}[theorem]{Lemma}

\theoremstyle{remark}
\newtheorem{remark}{Remark}

\theoremstyle{definition}

\title{\vspace{-3em}
A sharp Gaussian tail bound for sums of uniforms
}

\author{Xinjie He}

\author{Tomasz Tkocz}

\author{Katarzyna Wyczesany \\ \\ \emph{C\lowercase{arnegie} M\lowercase{ellon} U\lowercase{niversity}, P\lowercase{ittsburgh}, PA 15213, USA}}

\thanks{Research was supported in part by the NSF grant DMS-2246484. \\ Emails: \texttt{$\{$xinjieh, ttkocz, kwyczesa$\}$@andrew.cmu.edu}}

\begin{document}

\begin{abstract}
We prove that the tail probabilities of sums of independent uniform random variables, up to a multiplicative constant, are dominated by the Gaussian tail with matching variance and find the sharp constant for such stochastic domination.
\end{abstract}

\maketitle

\bigskip

\begin{footnotesize}

\end{footnotesize}

\bigskip

\section{Introduction}

Concentration inequalities have always played a prominent role in probability theory and beyond. They quantify the obvious intuition that random variables are likely to be close, or \emph{concentrate} around their mean. A prevalent source of concentration is independence, which has been studied to a great extent, see e.g. monographs \cite{BLM, Led}. Perhaps the most fundamental setting concerns sums of independent random variables, for which many robust tools have been developed. One of the precursors is undoubtedly Hoeffding's inequality from \cite{H}: for independent mean $0$ random variables $X_1, \dots, X_n$ with $|X_j| \leq 1$ a.s. for each $j$ and real numbers $a_1, \dots, a_n$ with $\sum_{j=1}^n a_j^2 = 1$, we have
\begin{equation}\label{eq:Hoeff}
\p{\left|\sum_{j=1}^n a_jX_j\right| > t} \leq 2e^{-t^2/2}, \qquad t > 0.
\end{equation}
In other words, a Gaussian tail bound of variance $1$ holds. When specialised to the Rademacher random variables (random signs), whose distribution is given by $\p{X_j = \pm 1} = \frac{1}{2}$, the sum $\sum a_jX_j$ has variance $1$. However, note that the tail of a standard Gaussian random variable (mean $0$, variance $1$) is asymptotic to $\frac{1}{\sqrt{2\pi}t}e^{-t^2/2}$, so the factor $\frac{1}{t}$ is \emph{missing}. Following Hoeffding's work, Efron remarked in \cite{Ef} that for applications to hypothesis testing, inequality \eqref{eq:Hoeff} ``is not sharp enough to be useful in practice'' and suggested that there is a universal constant $C$ such that for random signs we have
\begin{equation}\label{eq:efron}
\p{\left|\sum_{j=1}^n a_jX_j\right| > t} \leq C\p{|G| > t}, \qquad t > 0,
\end{equation}
where $G$ is a standard Gaussian random variable. This was proved by Pinelis in \cite{P1} and after some pursuit, the sharp constant has been found by Bentkus and Dzindzalieta in \cite{BDz}. Its value is approximately $3.18$ and comes from the case $n=2$, $a_1 = a_2 = \frac{1}{\sqrt{2}}$, $t = \sqrt{2}$. We refer, for instance, to Pinelis' work \cite{P3}, which besides asymptotic improvements (estimate \eqref{eq:efron} with $C = 1+O(1/t)$), gives  a detailed account on important milestones. For further multivariate analogues and extensions, we refer to \cite{CLT, NT, P6}. \emph{Truly} Gaussian tail bounds have also been established in a more general setting of martingales with bounded differences (i.e. bounds without \emph{missing factors}, as opposed to the classical Azuma inequality), see Pinelis' Theorem 5.4 in \cite{P5}, or Bentkus' work \cite{B1}.

This paper proves the sharp bound of the form \eqref{eq:efron} for sums of independent uniform random variables. Similarly to the Rademacher distribution playing a fundamental role in the class of all symmetric distributions, the uniform distribution is a building block of all (continuous) unimodal ones (every unimodal distribution is a mixture of uniform distributions, see, e.g. Chapter 1 in \cite{DJ}). To put it in a broader context, it is perhaps worth mentioning several earlier works devoted to various sharp bounds for sums of uniforms and their applications: \cite{Ball, Ball-wave, BK, BCh-epi, BChG, CET, CGT, CKT, ENT2, JT, LO}.

When the $X_j$ are chosen to be uniform on $[-1,1]$ with variance $\frac{1}{3}$, even the exponent of the bound by $e^{-t^2/2}$ in \eqref{eq:Hoeff} is suboptimal, as the Gaussian tail of matching variance would yield $e^{-3t^2/2}$. This can be rectified (see, e.g. Section 5 in \cite{BK}) but still lacks the missing factor $\frac1t$. On the other hand, the aforementioned works \cite{B1, P5} as well as Talagrand's paper \cite{Tal} retaining $\frac{1}{t}$ yield the main term $e^{-t^2/2}$, suboptimal for uniforms. Our main result is the following sharp bound, akin to Bentkus and Dzindzalieta's inequality from \cite{BDz}.

\begin{theorem}\label{thm:main}
Let $U_1, U_2, \dots$ be independent random variables uniform on $[-1,1]$. For every $n \geq 1$, real numbers $a_1, \dots, a_n$ with $\sum_{j=1}^n a_j^2 = 1$ and positive $t$, we have 
\begin{equation}\label{eq:main}
\p{\left|\sum_{j=1}^n a_jU_j\right| > t} \leq C_*\p{\frac{1}{\sqrt{3}}|G| > t},
\end{equation}
where $G$ is a standard Gaussian random variable (mean $0$, variance $1$) and the constant equals
\begin{equation}\label{eq:C*}
C_* = \sup_{0 < t < 1} \frac{1-t}{\p{|G|>t\sqrt{3}}} = 1.345118..
\end{equation}
(the supremum attained uniquely at $t_0 = 0.642908..$).
\end{theorem}

The value of the constant $C_*$ is clearly best possible because when $n=1$, $a_1 = 1$ and $t = t_0$, we have equality in \eqref{eq:main}. Plainly, the variance $\frac13$ of the Gaussian on the right hand side in \eqref{eq:main} matches the variance of the sum of uniforms $\sum a_jU_j$,  thus, in view of the central limit theorem, it is also optimal (i.e. cannot be replaced with any smaller one).

\section{Proofs}

\subsection{Overview}
In our proof of Theorem \ref{thm:main}, we use two different arguments, depending on whether $t$ is large or not. 
This strategy in \emph{low-resolution} mimics the one of Bentkus and Dzindzalieta from \cite{BDz}, however some fine points are different. In the small $t$ regime, \cite{BDz} relies on an improved Chebyshev's inequality (see Lemma 2.1 in \cite{BDz}), as well as Berry-Esseen bounds. For uniforms, this approach did not seem to yield satisfactory bounds, instead, we resort to log-concavity. More specifically, for the range $0 < t < 1$ which admits equality, we rely on arguments leveraging log-concavity, developed by Barthe and Koldobsky in \cite{BK}.

 For the range $t \geq 1$, we employ an inductive argument developed by Bobkov, G\"otze and Houdr\'e in \cite{BGH} for the Rademacher sums. This rests on a certain estimate for averages of the Gaussian tail function. We emphasise that the bound needed for the Gaussian tails in the case of uniforms is in fact stronger than the one for Rademachers (see Remark \ref{rem:ave-tail} in Section \ref{sec:rems}).

We now turn to the details for each of the two regimes.

\subsection{Approach for small $t$.}
Since the random vector $\tfrac12(U_1,\dots, U_n)$ is uniformly distributed on the unit volume cube $[-\tfrac12,\tfrac12]^n$ in $\R^n$, for a unit vector $a$ in $\R^n$ and $t > 0$,  we have
\[
\p{\bigg|\sum_{j=1}^n a_jU_j\bigg| \leq t} = \vol_{n} \bigg( \bigg\lbrace x \in [-\tfrac12,\tfrac12]^n, \ \bigg|\sum_{j=1}^n a_jx_j \bigg| \leq t/2\bigg\rbrace \bigg),
\]
that is, geometrically, this probability is equal to the volume (Lebesgue measure) of the section of the cube $[-\tfrac12,\tfrac12]^n$ by the slab $\{x \in \R^n,  \ |\sum a_jx_j| \leq t/2\}$ of width $t$. Barthe and Koldobsky have studied sections of minimal volume and have obtained the following sharp result (see Theorem 2 in \cite{BK}).

\begin{theorem}[Barthe-Koldobsky]\label{thm:BK}
Let $n \geq 1$. For every unit vector $a$ in $\R^n$ and $0 < t \leq \frac{3}{4}$,  we have
\[
\p{\bigg|\sum_{j=1}^n a_jU_j\bigg| \leq t}  \geq t.
\]
\end{theorem}

By the definition of $C_*$, this immediately gives \eqref{eq:main} for all $0 < t \leq \frac{3}{4}$.
We conjecture this to hold for all $t \leq 2(\sqrt{2}-1)$ as is known when $n=2$ (see Theorem 1 in \cite{BK}). Their arguments rest on the log-concavity of the sum of uniform random variables. Recall that a function $f\colon \R \to [0,+\infty)$ is called log-concave if it is of the form $f = e^{-\varphi}$ for a convex function $\varphi\colon \R \to (-\infty, +\infty]$ and a random variable is called log-concave if it has a density which is log-concave (see, e.g. Chapter 2.1 in \cite{BGVV} for background). Specifically, they relax the problem: ``given $n \geq 1$, $t > 0$, find $\min \p{\left|\sum_{j=1}^n a_jU_j\right| \leq t}$ over all unit vectors $a$ in $\R^n$'',
to the problem: ``given $t > 0$, find  $\min \p{|X| \leq t}$
 over all symmetric log-concave random variables $X$ with $\E X^2 = \frac13$.''
Using log-concavity, the latter problem naturally reduces to optimisation over the subfamily of truncated symmetric exponential densities. After explicit calculations, this leads to the following bound (obtained by combining Lemmas 4 and 5 from \cite{BK}). To state it, we define functions
\begin{align}
\psi(x) &= \frac{1}{x}\int_0^x \log^2(1-y) \dd y, \quad x > 0, \label{eq:def:psi}\\
G(t, p, x) &= 3t^2\psi(x) - \log^2(1-px), \quad t > 0, \ p, x \in (0,1). \label{eq:def:G}
\end{align}

\begin{lemma}[Barthe-Koldobsky]\label{lm:BK}
Let $t > 0$, $0 < p < 1$ and let $X$ be a symmetric log-concave random variable with $\E X^2 = \frac{1}{3}$. If
\[
\inf_{x \in (0,1)} G(t,p,x) \geq 0,
\]
then
\[
\p{|X| \leq t} \geq p.
\]
\end{lemma}

As this exact formulation, particularly suited for our purposes, is only implicit in their work, for completeness and reader's convenience, we sketch the argument in Section \ref{sec:lemmas}.

When applied to $X=\sum a_jU_j$, Theorem \ref{thm:BK} is then a consequence of the sharp bound 
\begin{equation}\label{eq:G-BK}
\inf_{x \in (0,1)} G(t,t,x) = 0,
\end{equation} 
which holds for all $ 0 < t \leq \frac34$ (Proposition 6 in \cite{BK}). However, this bound fails as soon as $t > \frac34$. Fortunately for us, the log-concave relaxation turns out to still be sufficient to handle our desired Gaussian bound \eqref{eq:main} for $\frac34 < t < 1$. In view of Lemma \ref{lm:BK}, to get \eqref{eq:main}, it is enough to prove the following technical lemma.

\begin{lemma}\label{lm:G}
For $t > 0$, set $p(t) = 1 - C_*\p{|G| > t\sqrt{3}}$, where $G$ is a standard Gaussian random variable. For every $\frac34 < t < 1$ and $0 < x < 1$, we have
\[
G(t,p(t),x) \geq 0.
\]
\end{lemma}

This inequality turns out to be quite subtle. We have not been able to find any easier solution than the rather non-elegant direct calculations leaning on convexity (and simple ``netting''). The proof is deferred to Section \ref{sec:lemmas}.

\subsection{Approach for large $t$.}
By the symmetry of the random variables involved, \eqref{eq:main} is equivalent to
\[
\p{\sum_{j=1}^n a_jU_j> t} \leq C_*\p{\frac{1}{\sqrt{3}}G > t}.
\] 
Using induction on $n$ and leveraging the independence of the $U_j$, we exploit the following natural approach going back to \cite{BGH}. We have,
\begin{equation}\label{eq:ind1}
\p{\sum_{j=1}^n a_jU_j > t}  = \p{\frac{\sum_{j=1}^{n-1} a_jU_j}{\sqrt{1-a_n^2}} > \frac{t - a_nU_n}{\sqrt{1-a_n^2}}}.
\end{equation}
Conditioning on the value of $U_n$, by the inductive hypothesis, as long as we have that $t - a_nU_n> 0$, we get
\begin{equation}\label{eq:ind2}
\p{\sum_{j=1}^n a_jU_j > t} \leq \E_{U_n}C_*\mathbb{P}_G\left(\frac{1}{\sqrt{3}}G > \frac{t - a_nU_n}{\sqrt{1-a_n^2}}\right).
\end{equation}
To finish this argument, it suffices to establish the following estimate on the averages of the Gaussian tail.

\begin{lemma}\label{lm:ave-tail}
Let $G$ be a standard Gaussian random variable. For every $0 < a < 1$ and $t > 1$, we have
\[
\frac{1}{2}\int_{-1}^1 \p{G > \frac{t+au}{\sqrt{1-a^2}}\sqrt{3}} \dd u \leq \p{G > t\sqrt{3}}.
\]
\end{lemma}

We postpone the proof to Section \ref{sec:lemmas}. 
Let us remark that the inductive step from \cite{BGH} for Rademacher random variables relied on the estimate
\begin{equation}\label{eq:BGH-ave-tail}
\frac{1}{2}\p{G > \frac{t\sqrt{3}+a}{\sqrt{1-a^2}}} + \frac{1}{2}\p{G > \frac{t\sqrt{3}-a}{\sqrt{1-a^2}}} \leq \p{G > t\sqrt{3}}.
\end{equation}
We will show in Remark \ref{rem:ave-tail} that Lemma \ref{lm:ave-tail} is stronger, in that the left hand side of \eqref{eq:BGH-ave-tail} is upper bounded by the left hand side from Lemma \ref{lm:ave-tail}.

\subsection{Proof of Theorem \ref{thm:main}}
We prove \eqref{eq:main} by induction on $n$. When $n=1$, $a_1 = \pm 1$ and consequently, the left hand side is simply $0$ for $t > 1$, so the inequality is trivial. For $0 < t \leq 1$, the inequality follows by the definition of $C_*$. 

For the inductive step, let $n \geq 2$ and suppose the result holds for every sum of $n-1$ uniform random variables. When $0 < t < 1$, combining Lemmas \ref{lm:BK} and \ref{lm:G} applied to $X = \sum_{j=1}^n a_jU_j$, which is symmetric log-concave with $\E X^2 = \frac13$, we get
\[
\p{|X| \leq t} \geq 1 - C_*\p{|G| > t\sqrt{3}}
\]
which is equivalent to \eqref{eq:main} and the argument is finished in this case. Thus, we can now assume that $t \geq 1$. Moreover, we can assume that $0 < a_n < 1$. Then $t - a_nU_n > 0$ a.s. and the inductive argument from \eqref{eq:ind1} and \eqref{eq:ind2} combined with Lemma \ref{lm:ave-tail} finishes the proof. \hfill$\square$

\section{Proofs of the auxiliary lemmas}\label{sec:lemmas}

\subsection{Proof of Lemma \ref{lm:BK} (sketch)}
Let $\mathcal{X}$ be the class of all symmetric log-concave random variables. First note that given arbitrary parameters $\sigma, t > 0$ and $0 < p < 1$, the following two statements are equivalent
\begin{align}\label{eq:state1}
&\text{For every $X \in \mathcal{X}$ with $\E X^2 = \sigma$, we have $\p{|X|\leq t} \geq p$}, \\\label{eq:state2}
&\text{For every $X \in \mathcal{X}$ with $\p{|X| \leq t} = p$, we have $\E X^2 \geq \sigma$.}
\end{align}
Indeed, if \eqref{eq:state1} holds and \eqref{eq:state2} does not, there is $X \in \mathcal{X}$ with $\p{|X| \leq t} = p$ and $\E X^2 < \sigma$. We consider $\tilde X =\lambda X$, $\lambda =  \sqrt{\frac{\sigma}{\E X^2}} > 1$. Then $\E \tilde X^2 = \sigma$, but $\p{|\tilde X| \leq t} = \p{|X| \leq t/\lambda} < \p{|X| \leq t} = p$, contradicting \eqref{eq:state1}. The converse is proved similarly. 

Using a standard argument of ``moving mass where it is beneficial'' for log-concave densities, Barthe and Koldobsky (see Lemma 4 in \cite{BK}) reduced the problem of finding $\inf \E X^2$ over all $X \in \mathcal{X}$ with $\p{|X| \leq t} = p$ to that over all $X$ with densities of the form $f(x) = ce^{-\alpha|x|}\1_{[-d,d]}(x)$ (truncated symmetric exponentials). After explicit calculations with such densities this is in turn reduced in their Lemma 5 to an optimisation problem over one parameter. Namely,
\[
\inf\{\E X^2, \ X \in \mathcal{X}, \p{|X| \leq t} = p\} = t^2\inf_{0 < x < 1} \frac{\psi(x)}{\log^2(1-px)},
\]
with $\psi$ defined in \eqref{eq:def:psi}. There is in fact a typo in their statement of Lemma 5: the square is missing at $\log(1-px)$). In view of the equivalence between \eqref{eq:state1} and \eqref{eq:state2}, showing that
\begin{equation}\label{eq:inf13}
t^2\inf_{0 < x < 1} \frac{\psi(x)}{\log^2(1-px)} \geq \frac{1}{3},
\end{equation}
results in $\p{|X| \leq t} \geq p$ for every $X \in \mathcal{X}$ with $\E X^2 = \frac{1}{3}$. \hfill$\square$

\subsection{Proof of Lemma \ref{lm:G}}
Our goal is to show that \eqref{eq:inf13} holds for all $\frac34 < t < 1$. That is, for all $0 < x < 1$ and $\frac34 < t < 1$ we have
\[
3t^2\psi(x) \geq \log^2(1-p(t)x),
\]
or, equivalently, after taking the square root on both sides,
\[
t\sqrt{3\psi(x)} \geq -\log(1-p(t)x).
\]
Note that the right hand side as a function of the quantity $p(t)$ is increasing. We first replace the function $p(t)$ by a piece-wise linear one which is pointwise larger for $\frac34 < t < 1$. With hindsight, we choose $t_0 = p_0 = \frac34$, $t_1 = 0.92$, $p_1 = 0.855$, $t_2 = 1$, $p_2 = 0.888$ and let
\[
\tilde p(t) = \begin{cases} p_0 + \frac{p_1-p_0}{t_1-t_0}(t-t_0), & t_0 \leq t \leq t_1, \\ p_1 + \frac{p_2-p_1}{t_2-t_1}(t-t_1), & t_1 \leq t \leq t_2. \\   \end{cases}
\]
This is the piecewise linear function which on $(t_j, t_{j+1})$ interpolates linearly between $(t_j, p_j)$ and $(t_{j+1}, p_{j+1})$, $j \in \{0, 1\}$.

\textbf{Claim.} \emph{We have, $p(t) < \tilde p(t)$ on $[\frac34, 1]$.}

\begin{proof}
Recall the definition, $p(t) = 1 - C_*\p{|G| > t\sqrt{3}}$. Since the density of $|G|$ is strictly decreasing, this is a concave function. Let $\ell_1(t), \ell_2(t)$ be its tangents at $t = 0.85$ and $t =1$ respectively, so that $p(t) \leq \ell_j(t)$, $j=1, 2$. We check that $\ell_1(t_0) < 0.748 < p_0$, $\ell_1(t_1) < 0.8545 < p_1$ and $\ell_2(t_1) < 0.8549 < p_1$, $\ell_2(t_2) < 0.888 = p_2$. By the construction of $\tilde p$, these bounds at the end-points finish the proof.   
\end{proof}

By the claim, we conclude that it suffices to show that for all $\frac34 < t < 1$ and $0 < x < 1$, we have
\[
t\sqrt{3\psi(x)} \geq -\log(1-\tilde p(t)x).
\]
Fix $0 < x < 1$. Observe that the right hand side as a function of $t$ is clearly convex on $[t_0, t_1]$ as well as on $[t_1, t_2]$. Since the left hand side is linear in $t$, once we have verified this inequality at $t = t_0, t_1, t_2$, we obtain that it is valid for all $t_0 \leq t \leq t_1$ and $t_1 \leq t \leq t_2$, as desired. Thus, our goal is to show that for $j = 0, 1, 2$ and all $0 < x  < 1$, we have
\begin{equation}\label{eq:goal-at-nodes}
t_j\sqrt{3\psi(x)} \geq -\log(1- p_jx).
\end{equation}
When $j=0$, $t_j = p_j = \frac34$, so \eqref{eq:goal-at-nodes} follows from \eqref{eq:G-BK} (Proposition 6 in \cite{BK}). For $j = 1, 2$, we begin by considering \emph{small} $x$. Recalling the definition of $\psi$ given in \eqref{eq:def:psi}, note that its Taylor series expansion at $x=0$ has the leading term $\frac{x^2}{3}$ and all coefficients positive (because $-\log(1-y)$ has positive coefficients, consequently, so does its square, $\log^2(1-y)$ and plainly, $\log^2(1-y) = y^2 + \dots$). As a result,
\begin{align}\label{eq:psi-convex}
&\psi(x) \ \text{is convex on (0,1)},\\
&\psi(x) > x^2/3, \quad 0 < x < 1.\label{eq:psi-lowbd}
\end{align}
Therefore, $\sqrt{3\psi(x)} > x$ and \eqref{eq:goal-at-nodes} is implied by $t_jx + \log(1-p_jx) \geq 0$. The left hand side as a function of $x$ is concave, it vanishes at $x = 0$ and we check that at $x = 0.15$, it equals $0.0007..$ and $0.007..$ for $j=1,2$, respectively. This proves \eqref{eq:goal-at-nodes} for all $0 < x < 0.15$. 

To finish the proof, we divide the interval $[0.15, 1]$ into $17$ smaller ones of equal length $0.05$, $[x_k, x_{k+1}]$, $x_k = 0.15 + 0.05k$, $k = 0, 1, \dots, 16$. We note that \eqref{eq:goal-at-nodes} is equivalent to
\begin{equation}\label{eq:goal-at-nodes-2}
\psi(x) \geq \frac{1}{3t_j^2}\log^2(1-p_jx).
\end{equation}
Let $g_j(x)$ be the right hand side. From \eqref{eq:psi-convex}, we know that $\psi(x)$ is convex, so is $g_j(x)$ (e.g. by the same argument as for $\psi$, its Taylor series expansion  at $x=0$ has all coefficients positive). On each interval $[x_k, x_{k+1}]$, we pointwise lower bound  $\psi(x)$ by its tangent $\tilde \psi_k(x)$ at the mid-point $\bar x_k = \frac{x_k+x_{k+1}}{2}$, that is $\tilde \psi_k(x) = \psi(\bar x_k) + \psi'(\bar x_k)(x-\bar x_k)$. By the convexity of $g_j(x)$, it remains to check that at the end-points $x = x_k, x_{k+1}$, we have $\tilde \psi_k(x) > g_j(x)$. Those calculations are gathered in Table \ref{tab} below.  \hfill$\square$
\begin{table}[H]
\begin{center}
\caption{Proof of \eqref{eq:goal-at-nodes-2}: lower bounds on the differences at the end-points between the linear approximations $\tilde \psi_k$ and $g_j$.}
\label{tab}
\hspace*{-7em}\begin{tabular}{r|cccccccccccccccccc}
$k$ & 0 & 1 & 2 & 3 & 4 & 5 & 6 & 7 & 8 & 9 & 10 & 11 & 12 & 13 & 14 & 15 & 16\\\hline
$10^3(\tilde \psi_k(x_k) - g_1(x_k))$ & $0.7$ & $1.4$ & $2.4$ & $3.6$ & $4.9$ & $6.4$ & $8.1$ & $9.8$ & $11$ & $13$ & $14$ & $16$ & $17$ & $20$ & $24$ & $34$ & $48$ \\ 
$10^3(\tilde \psi_k(x_{k+1}) - g_1(x_{k+1}))$ &  $1.5$ & $2.4$ & $3.6$ & $5$ & $6.5$ & $8.1$ & $9.9$ & $11$ & $13$ & $15$ & $16$ & $18$ & $21$ & $26$ & $41$ & $87$ & $304$ \\
$10^3(\tilde \psi_k(x_k) - g_2(x_k))$ & $1.3$ & $2.5$ & $4.2$ & $6.2$ & $8.5$ & $11$ & $14$ & $17$ & $20$ & $24$ & $27$ & $29$ & $31$ & $30$ & $28$ & $23$ & $3.0$ \\ 
$10^3(\tilde \psi_k(x_{k+1}) - g_2(x_{k+1}))$ &  $2.6$ & $4.2$ & $6.2$ & $8.6$ & $11$ & $14$ & $17$ & $20$ & $24$ & $27$ & $30$ & $31$ & $32$ & $30$ & $30$ & $41$ & $17$
\end{tabular}
\end{center}
\end{table}

\subsection{Proof of Lemma \ref{lm:ave-tail}}
Note that for $a = 0$, we have equality, hence it is enough to show that the left hand side as a function of $a$ is decreasing. Its derivative in $a$ equals
\[
-\frac{\sqrt{3}}{2(1-a^2)^{3/2}}\int_{-1}^1 \phi\left( \frac{t+au}{\sqrt{1-a^2}}\sqrt{3} \right)(u+at)\dd u
\]
where $\phi(x) = \frac{1}{\sqrt{2\pi}}e^{-x^2/2}$ is the density of $G$. We will show that this is negative for all $t > 1$ and $0 < a < 1$. This is clear when $at \geq 1$, since then $u + at \geq 0$ for all $-1 \leq u \leq 1$. From now on, we assume that $at < 1$. 
Let $t_\pm = \frac{t\pm a}{\sqrt{1-a^2}}\sqrt{3}$. After a change of variables, $v = \frac{t + au}{\sqrt{1-a^2}}\sqrt{3}$ and using that $v\phi(v) = -\phi'(v)$, it remains to show that
\[
h(a,t) = \frac{\phi(t_-)-\phi(t_+)}{t\sqrt{3}\sqrt{1-a^2}} - \int_{t_-}^{t_+} \phi(v) \dd v > 0, \qquad 0 < a < 1, \ 1 < t < \frac{1}{a}.
\]
Again, there is equality at $a=0$, so it is enough to show that the derivative in $a$ of $h$ is positive. We have,
\[
\frac{\partial h}{\partial a} = \frac{a\phi(t_-)}{t\sqrt{3}(1-a^2)^{5/2}}\left[(3a^2t^2+6at+a^2+2)e^{-\frac{6at}{1-a^2}}-(3a^2t^2-6at+a^2+2)\right].
\]
We let $b = at$ and claim that for every $0 < a < b  < 1$,
\[
\exp\left(-\frac{6b}{1-a^2}\right) > \frac{3b^2-6b+a^2+2}{3b^2+6b+a^2+2} = 1 - \frac{12b}{3b^2+6b+a^2+2},
\]
which will finish the proof. Notice that the left hand side is clearly decreasing in $a$, whilst the right side side is increasing, thus it suffices to show this inequality for $a=b$, that is
\[
\exp\left(-\frac{6a}{1-a^2}\right) > \frac{2a^2-3a+1}{2a^2+3a+1} = \frac{(1-2a)(1-a)}{(1+2a)(1+a)},\qquad 0 < a < 1.
\]
When $a \geq \frac12$, the right hand side is nonpositve, so the inequality is obvious. For $0 < a < \frac12$, equivalently, we would like to argue that
\[
f(a) = -\frac{6a}{1-a^2} - \log\left(\frac{(1-2a)(1-a)}{(1+2a)(1+a)} \right)
\]
is positive. We have $f(0) = 0$ and
\[
f'(a) = \frac{36a^4}{(1-a^2)^2(1-4a^2)},
\]
which is clearly positive, completing the argument. \hfill$\square$

\section{Final remarks}\label{sec:rems}

\begin{remark}\label{rem:suboptimal-C}
Our main result, inequality \eqref{eq:main} with a suboptimal value of constant $C_*$ can be obtained from a majorisation result of Theorem 2 from \cite{ENT2}, combined with Pinelis' technique ``from (generalised) moments to tails'', \cite{P3}. Such an argument yields \eqref{eq:main} with $c = e^2/2 = 3.69..$ in place of $C_*$. To sketch this argument, note first that for every even function $f$, nondecreasing on $[0, +\infty)$, we have by Markov's inequality,
\[
\p{\big|\sum a_jU_j\big| > t} \leq \frac{1}{f(t)}\E f\big(\sum a_jU_j\big).
\]
If $f$ is additionally $C^1$ with $f'$ convex on $[0,+\infty)$, the sum of uniforms can be dominated by the sum of Gaussians with matching variances (Theorem 6 combined with Remark 18 from \cite{ENT2}), resulting with the bound
\[
\p{\big|\sum a_jU_j\big| > t} \leq \frac{1}{f(t)}\E f(G/\sqrt{3}).
\]
Letting $f_u(x) = (|x|-u)_+$ with a carefully chosen $u$ (depending on $t$) allows to estimate this quantity by the tail function (leveraging only the log-concavity of $G$),
\[
\frac{1}{f_u(t)}\E f_u(G/\sqrt{3}) \leq \frac{e^2}{2}\p{\frac{1}{\sqrt{3}}|G| > t}
\]
(see Theorem 3.11 in \cite{P5}, or the direct computation from \cite{CLT} following Claim 2).
\end{remark}

\begin{remark}\label{rem:ave-tail}
Let $G, U, \varepsilon$ be independent, $G$ standard Gaussian, $U$ uniform on $[-1,1]$ and $\varepsilon$ uniform on $\{-1,1\}$. It can be checked that the following function
\[
x \mapsto \p{G > \frac{t\sqrt{3} + a\varepsilon\sqrt{3x}}{\sqrt{1-a^2}}}
\]
is convex on $[0,1]$. As a result, by Jensen's inequality,
\[
\p{G > \frac{t\sqrt{3} + a\varepsilon}{\sqrt{1-a^2}}} = \p{G > \frac{t\sqrt{3} + a\varepsilon\sqrt{3\E U^2}}{\sqrt{1-a^2}}} \leq \E_U\p{G > \frac{t\sqrt{3} + a\varepsilon|U|\sqrt{3}}{\sqrt{1-a^2}}}.
\]
Since $\varepsilon|U|$ has the same distribution as $U$, we conclude that Lemma \ref{lm:ave-tail} implies \eqref{eq:BGH-ave-tail}.
\end{remark}

\begin{remark}\label{rem:self-norm}
As \eqref{eq:efron} has direct applications to tail bounds for self-normalising sums of \emph{symmetric} random variables (see \cite{Ef, Pena}), our inequality \eqref{eq:main} finds a similar use for sums of symmetric \emph{unimodal} ones. Let $X_1, \dots, X_n$ be independent symmetric unimodal random variables, that is, each $X_i$ has the same distribution as $R_jU_j$, where $R_1, \dots, R_n$ are some independent nonnegative random variables, independent of $U_1, \dots, U_n$ which are i.i.d. uniform on $[-1,1]$ (we refer to \cite{DJ} for basic background, and in particular, Theorem 1.5 for such representation). Then, by conditioning on the values of the $R_j$, a direct consequence of Theorem \ref{thm:main} is the following bound
\[
\p{\left|\frac{\sum_{j=1}^n R_jU_j}{\sqrt{\sum_{j=1}^n R_j^2}}\right| > t} \leq C_*\p{\frac{1}{\sqrt{3}}|G|>t}, \qquad t > 0.
\]
We emphasise that when the $X_1, X_2, \dots$ are i.i.d. with finite variance, the self-normalised sum 
\[
\frac{\sum_{j=1}^n R_jU_j}{\sqrt{\sum_{j=1}^n R_j^2}} = 
\frac{ \frac{\sum_{j=1}^n R_jU_j}{\sqrt{n}} }{ \sqrt{\frac{\sum_{j=1}^n R_j^2}{n}} }
\]
converges in distribution to $\frac{G\sqrt{\E R_1^2U_1^2}}{\sqrt{\E R_1^2}} = \frac{1}{\sqrt{3}}G$.
\end{remark}

\begin{remark}\label{rem:complex}
K\"onig and Koldobsky in \cite{KK2} proved an analogue of Theorem \ref{thm:BK} in the complex setting where the $U_j$ are uniform on the unit disc in the complex plane. Thus it will be of interest to find the best constant for higher dimensional analogues of \eqref{eq:main}, where the $U_j$ are uniform on Euclidean balls or spheres. Such tail comparisons are known with (subobtimal) \emph{universal} constants, which have been established independently in \cite{NT, P6}. 
\end{remark}

\begin{remark}\label{rem:BK-short-proof}
There is a short argument justifying Theorem \ref{thm:BK} for a slightly smaller range: $0 < t \leq \frac{2}{3}$ (still containing the value of $t = t_0 = 0.64..$ where \eqref{eq:main} is tight). It relies on a curious formula with the negative moment, namely
\[
\p{\bigg|\sum_{j=1}^n a_jU_j\bigg| \leq t} = t\E\left[\bigg| t\xi_0 + \sum_{j=1}^n a_j\xi_j\bigg|^{-1}\right], \qquad t > 0,
\]
holds for every unit vector $a \in \R^n$, where $\xi_0, \xi_1, \dots$ are i.i.d. random vectors each uniform on the unit Euclidean sphere in $\R^3$ and $|\cdot|$ denotes the Euclidean norm. This formula goes back to \cite{KK}, Proposition 2, where it was derived using the Fourier transform (a different argument rests on writing densities at $0$ using negative moments, see e.g. (1) in \cite{CKT}, or \cite{CNT}). Let $X = \sum a_j\xi_j$. By rotational invariance and independence,
\[
\E\left[\left|t\xi_0 + X\right|^{-1}\right] = \E\min\left\{t^{-1},\left|X\right|^{-1}\right\} = \E f_t(|X|^2),
\]
(see, e.g. Lemma 2 in \cite{GTW} for a more general statement in arbitrary dimension), where we set
\[
f_t(x) = \min\{t^{-1},x^{-1/2}\}.
\]
As long as $0 < t \leq \frac{2}{3}$, for the tangent to $x^{-1/2}$ at $x = 1$, $g(x) = 1 - \tfrac12(x-1)$ we have, $g(x) \leq f_t(x)$, $x \geq 0$ (by convexity, it suffices to check $g(0) \leq f_t(0)$ which is $g(0) = \frac{3}{2} \leq \frac1t= f_t(0)$). Since $\E |X|^2 = 1$, we get
\[
\E f_t(|X|^2) \geq \E g(|X|^2) = 1 - \frac12(\E |X|^2 - 1) = 1,
\]
which, in view of the previous formula, gives Theorem \ref{thm:BK} for all $0 < t \leq \frac23$.  We finish with recalling V. Milman's beautiful conjecture that given $t > 0$, the minimum of $\p{|\sum a_j U_j| \leq t}$ over all unit vectors $a$ in $\R^n$ is attained at $a$ of the form $(1, \dots, 1, 0, \dots, 0)$ for a suitable number of $1$s (see \cite{BK}, \cite{KK} and \cite{MSZZ} for progress made so far). In light of this, we conjecture that for all unit vectors $a$ in $\R^n$,
\[
\p{\left|\sum_{j=1}^n a_jU_j\right| > t} \leq \p{\frac{1}{\sqrt{3}}|G| > t}, \qquad t > t_1 = 1.04288..,
\]
where $t_1$ is the unique solution of $\p{|\frac{U_1+U_2+U_3}{\sqrt{3}}| > t} = \p{\frac{1}{\sqrt{3}}|G| > t}$, $t > 0$. There are other natural probabilistic-geometric questions featuring similar  intermediate dimensional symmetry breaking, see e.g. \cite{HJN} for the so-called propeller conjecture, \cite{BN} for Khinichin-type inequalities and links to projections, or \cite{CNT} for links to sections, in particular Conjecture 2  in \cite{NT-surv}. 
\end{remark}

\end{document}